\documentclass[11pt]{article}
\usepackage{amssymb}
\usepackage{epsfig}

\newtheorem{Th}{Theorem}[section]

\def\be{\begin{equation}}
\def\ee{\end{equation}}
\def\bea{\begin{eqnarray}}
\def\eea{\end{eqnarray}}
\def\bes{\begin{eqnarray*}}
\def\ees{\end{eqnarray*}}

\def\nn{\nonumber}

\def\<{\langle}
\def\>{\rangle}
\def\lb{\label}
\def\bs{\setminus}

\def\R{{\bf R}}

\def\Z{{\bf Z}}

\def\N{{\bf N}}
\def\U{{\bf U}}

\def\Q{{\bf Q}}

\def\bb{{\beta}}
\def\ga{{\gamma}}

\def\th{{\theta}}

\def\Om{{\Omega}}
\def\ep{{\epsilon}}
\def\lm{{\lambda}}

\def\Dl{{\Delta}}
\def\dl{{\delta}}
\def\sg{{\sigma}}
\def\dm{{\diamond}}

\def\P{{\cal P}}

\def\Sp{{\rm Sp}}

\def\hb{\vrule height0.18cm width0.14cm $\,$}

\title{On the common index jump theorem and further developments}

\author{Huagui Duan$^{1, 2}$ \thanks{Partially supported by National Key R\&D Program of China (2020YFA0713300), NSFC (12271268 and 12361141812) and Natural Science Foundation of Tianjin (25JCZDJC01030).  E-mail: duanhg@nankai.edu.cn.},
\  Yiming Long$^{3}$ \thanks{Partially supported by National Key R\&D Program of China (2020YFA0713300), NSFC
(11131004, 11671215 and 11790271), Nankai University, Wenzhong
Foundation and Nankai Zhide Foundation. E-mail: longym@nankai.edu.cn.},
\  Wei Wang$^{4}$ \thanks{Partially supported by NSFC (12025101). E-mail: wangwei@math.pku.edu.cn},
\  Chaofeng Zhu$^{5}$ \thanks{Partially supported by National Key R\&D Program of China (2020YFA0713300),
NSFC (12271268), Nankai Zhide Foundation and Nankai University.}\\\\
$^{1}$ Department of Mathematics, Jilin University, Changchun, Jilin, 130012\\
$^{2}$ School of Mathematical Sciences and LPMC, Nankai University, Tianjin 300071\\
$^{3, 5}$ Chern Institute of Mathematics and LPMC, Nankai University, Tianjin 300071\\
$^{4}$ School of Mathematical Sciences, Peking University, Beijing 100871\\
The People's Republic of China\\}

\begin{document}
\date{}
\maketitle

\begin{abstract}
In \cite{LZ02} published in ``Annals of Math.'', Long and Zhu established originally the common index jump theorem
(CIJT) for symplectic paths in 2002, which has played an important role in later studies on periodic solution orbits
for Hamiltonian systems, Reeb flows, and geodesic problems. This (CIJT) was generalized to its enhanced version
(ECIJT) by Duan, Long and Wang in \cite{DLW16} in 2016. Started from \cite{GG20} of 2020, and finally in \cite{CGG24}
of 2024, a similar index theorem was obtained, i.e., Theorem 3.3 of \cite{CGG24}, which was given the name ``index
recurrence theorem'' there. In this short note, we give detailed proofs to show that the major assertions, i.e., the
first 4 assertions in the total of 5 assertions, in Theorem 3.3 of \cite{CGG24} as well as all the assertions in
\cite{GG20} actually coincide completely with results in (ECIJT) of \cite{DLW16}.
\end{abstract}



Abbreviated title: Common index jump theorem and its extensions.

\renewcommand{\theequation}{\thesection.\arabic{equation}}
\renewcommand{\thefigure}{\thesection.\arabic{figure}}

\setcounter{section}{0}
\setcounter{equation}{0}\setcounter{figure}{0}
\section{Introduction} %Section 1

In the paper \cite{LZ02} of Long and Zhu published in 2002, the common index jump theorem (CIJT for short below)
was established, and since then this theorem has become a powerful tool in the study of many problems related to
periodic solutions of Hamiltonian systems.

Note that in the paper \cite{DLW16} of Duan, Long and Wang published in ``Calculus of Variations" in 2016, CIJT
was improved to the enhanced common index jump theorem (ECIJT for short below), i.e., Theorem 3.5 of \cite{DLW16},
with some more precise formula on index iterations. Then in recent \cite{DLLW24} of Duan, Liu, Long and Wang in
2024, ECIJT was further generalized.

On the other hand, In \cite{GGM18} of Ginzburg, G\"urel and Macarini published in 2018, a so-called index
recurrence theorem, i.e., Theorem 4.1 there, was obtained, and the authors wrote in their abstract of \cite{GGM18}
specially that ``{\it On the combinatorial side of the question, we revisit and reprove the enhanced common jump
theorem (i.e., the (ECIJT) Theorem 3.5) from \cite{DLW16} and interpret it as an index recurrence result (i.e.,
the Theorem 4.1 in \cite{GGM18})}".

In \cite{GG20} of Ginzburg and G\"urel published in 2020, its Theorem 5.2  was obtained, which was also called
the index recurrence theorem. Note that although the Theorem 4.1 of \cite{GGM18} assumed the strongly
non-degenerate condition on the symplectic paths studied, but the Theorem 3.5 and Remark 3.6 in \cite{DLW16}
published in 2016 used no any non-degeneracy conditions and allow degenerate symplectic paths, which thus covered
already completely all the assertions in Theorem 5.2 of \cite{GG20} published in 2020, which is 4 years later than
\cite{DLW16}.

In \cite{CGG24} of Cineli, Ginzburg and G\"urel appeared as arXiv:2410.13093v3 on Oct. 11, 2025, the Theorem 5.2 of \cite{GG20} was extended to the Theorem 3.3 of \cite{CGG24}
and was also called index recurrence theorem too, whose proof was a replenished version of and was more rigorous
than that of Theorem 5.2 in \cite{GG20}. This Theorem 3.3 became one of the main tools in the proof of the main
multiplicity result on closed Reeb flow orbits in \cite{CGG24}.

Note that the Theorem 3.3 in \cite{CGG24} contains 5 assertions (IR1)-(IR5) in total, and it needs to be pointed
out that the first 4 assertions there coincide completely with the identities (3.20)-(3.22) of Theorem 3.5
(ECIJT) and estimates (3.42) and (3.43) in Remark 3.6 in \cite{DLW16} already published in 2016. Note also that
the proof of (IR1)-(IR4) in \cite{CGG24} occupied the major part in the proofs of Theorem 3.3 in \cite{CGG24}.
On the other hand because all the assertions in Theorem 5.2 in \cite{GG20} are covered by these (IR1)-(IR4) in
\cite{CGG24}, as mentioned above all the assertions of Theorem 5.2 in \cite{GG20} coincide also with those
proved in \cite{DLW16} published in 2016 already too.

Because for notations and proofs, in \cite{CGG24} the authors used the Arnold-Givental normal forms of real
quadratic forms established by Arnold and Givental in \cite{AG} of 2001, and in \cite{DLW16} the authors used
the normal forms and basic normal forms of real symplectic matrices established by Dong, Han and Long in
\cite{LD00}, \cite{LH99} and \cite{Lo00} (cf. Subsections 1.4-1.9 on pp.16-47 and Chapter 8 in \cite{Lo02} for
details), the coincidences between above mentioned (IR1)-(IR4) of Theorem 3.3 in \cite{CGG24} and results
already proved in \cite{DLW16} are not that obvious for usual readers. This short paper is devoted to give
detailed proofs for such coincidences, which are contained specially in the 4 Propositions below.

In the following, when the page numbers of \cite{CGG24} appear we are using those in arXiv:2410.13093v3 on Oct. 11, 2025.

\medskip

\setcounter{equation}{0}\setcounter{figure}{0}
\section{Theorem 3.3 and Remark 3.6 of \cite{DLW16} as well as Theorem 3.3 of \cite{CGG24}} %Section 2

Note that in Theorem 3.5 on p.145 of \cite{DLW16}, the authors proved the following theorem,

{\bf Theorem 3.5 in \cite{DLW16}.} (The enhanced common index jump theorem for symplectic paths) {\it Let
$\ga_k\in \P_{\tau_k}(2n)$ for $k=1,\ldots,q$ be a finite collection of symplectic paths. Let
$M_k = \ga_k(\tau_k)$. We extend $\ga_k$ to $[0,+\infty)$ by (3.9) of \cite{DLW16} inductively. Suppose
$$    \hat{i}(\ga_k,1) > 0, \qquad \forall\;k = 1, \ldots, q. \hskip 8 cm {\rm (3.17)\;of\;[3]}. $$
Let
$$   \check{m} \equiv \check{m}(\ga_1, \ldots, \ga_q) = \min\{\check{M}_k\;|\;1\le k\le q\}. \hskip 6 cm {\rm (3.18)\;of\;[3]}. $$
where
$$ \check{M}_k = \left\{\matrix{
      \min\{m\in\N\;|\;m\th\in 2\pi\N,\; e^{\sqrt{-1}\th}\in\sg(M_k)\;{\rm with}\;\th\in (0,2\pi)\cap \pi \Q\},
                & \;{\rm if}\;M_k\in \Sp_{vnu}(2n), \cr
      + \infty, & \;{\rm if}\;M_k\in \Sp_{cnu}(2n), \cr}\right.  $$
and
\bea
\Sp_{cnu}(2n) &=& \{M\in\Sp(2n)\;|\;\dim(M^m-I) = \dim(M-I),\;\forall\;m\in\N\}, \nn\\
\Sp_{vnu}(2n) &=& \Sp(2n)\bs \Sp_{cnu}(2n). \nn\eea

Then for every integer $\bar{m}\in\N$, there exist infinitely many $(q+1)$-tuples $(N,m_1,\ldots, m_q)\in \N^{q+1}$
such that
$$ \nu(\ga_k,2m_k-m) = \nu(\ga_k,2m_k+m) = \nu(\ga_k,1),\;\forall\;1\le k\le q,\;1\le m < \check{m}, \hskip 1 cm {\rm (3.19)\;of\;[3]}, $$
and the following hold for all $1\le k\le q$ and $1\le m\le \bar{m}$,
\bea
&& \nu(\ga_k,2m_k-m) = \nu(\ga_k,2m_k+m) = \nu(\ga_k,m),  \hskip 4.5 cm {\rm (3.20)\;of\;[3]}, \nn\\
&& i(\ga_k,2m_k+m) = 2N + i(\ga_k,m),  \hskip 6.6 cm {\rm (3.21)\;of\;[3]}, \nn\\
&& i(\ga_k,2m_k-m) = 2N - i(\ga_k,m) - 2(S_{M_k}^+(1)+Q_k(m)), \hskip 3 cm {\rm (3.22)\;of\;[3]}, \nn\\
&& i(\ga_k,2m_k) = 2N - (S_{M_k}^+(1)+C(M_k)-2\Dl_k), \hskip 4.5 cm {\rm (3.23)\;of\;[3]}, \nn\eea
where as in \cite{LZ02}, we let
$$  \Dl_k = \sum_{0<\{m_k\th/\pi\}<\dl}S_{M_k}^-(e^{\sqrt{-1}\th}),  \hskip 7.5 cm {\rm (3.24)\;of\;[3]}, $$
and we define }
$$  Q_k(m) = \sum_{\th\in (0,2\pi),\;e^{\sqrt{-1}\th}\in \sg(M_k)
         \atop \{\frac{m_k\th}{\pi}\} = \{\frac{m\th}{2\pi}\} = 0}S_{M_k}^-(e^{\sqrt{-1}\th}).
            \qquad \hskip 5 cm {\rm (3.25)\;of\;[3]}$$

Note that in (3.42) and (3.43) contained in Remark 3.6 on p.145 in \cite{DLW16}, the authors proved also the
following two estimates.

{\bf Remark 3.6 (ii) in \cite{DLW16}.} {\it By (4.10) and (4.40) in \cite{LZ02} (cf. (11.1.10) and (11.2.14)
of \cite{Lo02}), we have
$$   m_k = \left(\left[\frac{N}{\bar{M}\hat{i}(\ga_k,1)}\right]+\chi_k\right)\bar{M},
                  \qquad 1\le k\le q, \hskip 3 cm {\rm (3.42)\;of\;[3]} $$
where $\chi_k=0$ or $1$ for $1\le k\le q$, $\bar{M}$ is a positive integer such that $\frac{\bar{M}\th}{\pi} \in\Z$
whenever $e^{\sqrt{-1}\th}\in \sg(M_k)$ and $\frac{\th}{\pi}\in \Q$ for some $1\le k\le q$, and we set $\bar{M}=1$
if no such eigenvalues exist. Furthermore, by (4.20) in Theorem 4.1 of \cite{LZ02} (cf.(11.1.20) of \cite{Lo02}),
for any $\ep>0$, we can choose $N$ and $\{\chi_k\}_{1\le k\le q}$ such that }
$$   \left|\left\{\frac{N}{\bar{M}\hat{i}(\ga_k,1)}\right\} - \chi_k\right| < \ep, \qquad 1\le k\le q.
           \hskip 4 cm {\rm (3.43)\;of\;[3]} $$

In the above, $[a]=\max\{k\in\Z\;|\;k\le a\}$ and $\{a\}=a-[a]$ for all $a\in\R$.

\medskip

On the other hand, in p.31 of \cite{CGG24}, under similar
conditions the authors proved their Theorem 3.3 which contains five assertions (IR1)-(IR5). Here we are interested
in the following first four assertions (IR1)-(IR4) on the index iteration theories of symplectic matrix paths started
from the identity matrix.

{\bf Theorem 3.3 in \cite{CGG24}.} {\it For any $\eta>0$ and any $\ell_0\in\N$, there exists a sequence $C_l\to \infty$
of real numbers and integer sequences $d_{il}\to \infty$ and $k_{ijl}\to \infty$ (as $l\to\infty$) such that for
all $i$, $j$ and $l$, and all $\ell\in\Z$ in the range $1\le |\ell|\le \ell_0$, we have
\bea
&& (IR1) \qquad |\hat{\mu}(\Phi_{ij}^{k_{ijl}}) - d_{il}| < \eta,   \hskip 9 cm (IR1-1)\nn\\
&& \qquad\qquad {\it and\;in\;particular,}\;d_{il} = [\hat{\mu}(\Phi_{ij}^{k_{ijl}})],\;{\it and}\; \nn\\
&& \qquad\qquad d_{il} - m \le \mu_-(\Phi_{ij}^{k_{ijl}}) \le \mu_+(\Phi_{ij}^{k_{ijl}}) \le d_{il}+ m; \hskip 5.2 cm (IR1-2)\nn\\
&& (IR2) \qquad \mu_{\pm}(\Phi_{ij}^{k_{ijl}+\ell}) = d_{il} + \mu_{\pm}(\Phi_{ij}^{\ell}) \qquad {\it when}\quad 0<\ell\le \ell_0; \nn\\
&& (IR3) \qquad \mu_+(\Phi_{ij}^{k_{ijl}-\ell}) = d_{il} - \mu_-(\Phi_{ij}^{\ell}) + (\bb_+(\Phi_{ij}^{\ell}) - \bb_-(\Phi_{ij}^{\ell})), \nn\\
&& \qquad\qquad {\it where}\;0<\ell\le \ell_0<k_{ijl}\; {\it and}\;\bb_+ - \bb_- \le m\; {\it with}\;\bb_{\pm}=0\;{\it when}\;
                  \Phi^{\ell}_{ij}\;{\it is~ non-degenerate};   \nn\\
&& (IR4) \qquad {\it Assertions\;(IR1)-(IR3)\;continue\;to\;hold\;with\;the\;same\;values\;of\;}d_{il}\;
         {\it and}\;k_{ijl}\;{\it when} \nn\\
&& \qquad\qquad \Phi_{ij}\;{\it is\;replaced\;by\;its\;non-degenerate\;part\;}
         \Psi\;{\it with\;all\;invariants\;of\;}\Phi_{ij}\;{\it and\;its\;}   \nn\\
&& \qquad\qquad {\it iterates\;replaced\;by\;their\;counterparts\;for\;}\Psi,\;{\it i.e.,}\;m\;{\it replaced\;by\;}
           m',\;\mu_{\pm}(\Phi_{ij}^{k_{ijl}})\; \nn\\
&& \qquad\qquad {\it replaced\;by\;}\mu(\Phi_{ij}^{k_{ijl}}),\;\bb_{\pm}\;{\it set\;to\;be\;}0,\;{\it etc.},  \nn\\
&& (IR5) \qquad C_l-\eta < k_{ijl}a_{ij} < C_l, \;{\it and}\; ka_{ij} < C_l-\eta\; {\it when}\;k<k_{ijl}\;{\it and}\;
                ka_{ijl} > C_l\; {\it when}\; k>k_{ijl}.  \nn\eea}

Note that all of (IR1)-(IR4) are in fact claims on only one single symplectic path, and clearly (IR4) is a simple
consequence of (IR2)-(IR3). Thus the authors of \cite{CGG24} pointed out specially that their proofs for (IR1)-(IR3)
can be simplified to the following three properties on a single symplectic path.
\bea
(IR1") && |\hat{\mu}(\Phi^{k}) - d| < \eta, \qquad{\it and}\;\;   \hskip 7.8 cm (IR1"-1)\nn\\
       &&  d - m \le \mu_-(\Phi^{k}) \le \mu_+(\Phi^{k}) \le d + m; \hskip 6 cm (IR1"-2)\nn\\
(IR2") && \mu_{\pm}(\Phi^{k+\ell}) = d + \mu_{\pm}(\Phi^{\ell}) \qquad {\it when}\quad 0<\ell\le \ell_0, \nn\\
(IR3") && \mu_+(\Phi^{k-\ell}) = d - \mu_-(\Phi^{\ell}) + (\bb_+(\Phi^{\ell}) - \bb_-(\Phi^{\ell}), \nn\eea
{\it where $0<\ell\le \ell_0<k$ and $\bb_+-\bb_- = b_+-b_- \le m$ with $\bb_{\pm}=0$ when $\Phi^\ell$ is non-degenerate. }

\medskip

Next we explain the details of the comparisons of results contained in the papers \cite{CGG24} and \cite{DLW16}
related to (IR1)-(IR4) in the following 4 sections.

\medskip\medskip

\section{The coincidence of $\bb_{-}(\ga)$ in \cite{CGG24} and $S_{\ga(1)}^-(1)$ in Long's book
\cite{Lo02} of 2002} %Section 3

To get precise definitions of certain quantities used in \cite{CGG24} and understand the proof of Theorem 3.3
in \cite{CGG24}, we make first the following concrete computations, and then prove the Proposition 1 below.

As in p.28 of \cite{CGG24}, using the normal form of symplectic matrices in Section 2.4 of \cite{AG}, the
authors considered a totally degenerate matrix $A\in \Sp(2m)$. Then $A$ can be written into $A=\exp(JQ)$,
where all eigenvalues of $JQ$ are zero. The quadratic form $Q$ can be symplectically decomposed into a sum
of terms of the following three types according to Section 2.4 of \cite{AG}:

(i) the identically zero quadratic form on $\R^{2\nu_0}$,

(ii) the quadratic form $Q_0 = p_1q_2 + p_2q_3 + \cdots + p_{d-1}q_d$ in Darboux coordinates on $\R^{2d}$,
where $d\ge 1$ is odd,

(iii) the quadratic forms $Q_{\pm} = \pm (Q_0+p_d^2/2)$ on $\R^{2d}$ for any $d\ge 1$.

Note that when the integer $d$ in the above (ii) is even, the quadratic form $Q_0$ can be split into a sum of
sub-quadratic forms $Q_i$ defined on $\R^{2d_i}$ for $1\le i\le i_0$ for some integer $i_0\ge 2$ with each
$d_i$ is an odd integer.

\medskip

{\bf Step 1.1.} {\it On the invariants $b_{\ast}(Q)$ and $\bb_{\pm}(Q)$.}

In \cite{CGG24}, for $\ast = 0, \pm$, the authors defined
\be   b_{\ast}(Q) = \{ \rm the\;number\;of\;Q_{\ast}\;{\rm appeared\;in\;the\;decomposition\;of\;}Q\}.  \lb{3.1}\ee
Note that here $\nu_0(Q)$ is defined to be the number of identically zero quadratic forms in the decomposition
of $Q$ on $\R^{2\nu_0}$.

Therefore for $A=\exp(JQ)$ there holds
\be   \nu(A) = \dim\ker (A-I) = 2(b_0 +\nu_0) + b_+ + b_-,  \lb{3.2}\ee
where $\nu(A)$ is the geometric multiplicity of the eigenvalue $1$ of $A$ as usual (cf. Definition 5.1.1 on
p.111 of \cite{Lo02}).

As shown below in our discussion, such a quadratic form $Q$ corresponds to a symplectic path
$\Phi\in C([0,1],\Sp(2n))$ up to a homotopy. Then as described in p.28 of \cite{CGG24}, the authors decompose $\Phi$ into a sum of two paths $\Phi = \Phi_0\oplus \Psi$ such that $\Phi_0(1)\in \Sp(2n_0)$ is totally degenerate
and $\Psi(1)\in \Sp(2n_1)$ is non-degenerate. In particular $n_0$ is the algebraic multiplicity of the
eigenvalue $1$ of $\Phi(1)$ and $n_0+n_1=n$. Then the authors of \cite{CGG24} defined first
$$ b_{\ast}(\Phi) = b_{\ast}(\Phi_0(1))\quad {\rm for}\;\;\ast = 0,\pm \qquad
                     {\rm and}\quad \nu_0(\Phi) = \nu_0(\Phi_0(1)). $$
Then they defined
\be   \bb_{\pm}(\Phi) = \nu_0(\Phi) + b_0(\Phi) + b_{\pm}(\Phi).  \lb{3.3}\ee

\medskip

{\bf Step 1.2.} {\it On $b_{0}(\Phi)$.}

We study the quadratic form $Q_{0}$ on $\R^{2d}$ with an odd integer $d\ge 1$ first.

Note that the values of the function $b_0$ introduced in p.28 of \cite{CGG24} is defined by $b_0(Q)=\nu_0(Q_0)$.
From the definition of $Q_0$ in the above (ii), we obtain
\be  Q_0 = \left(\matrix{ O & A\cr
                         A^T & O\cr}\right)_{2d\times 2d}, \quad{\rm with\;}
 A = \left(\matrix{0 & 1 & 0 & \cdots & 0 & 0\cr
                   0 & 0 & 1 & \cdots & 0 & 0\cr
                   0 & 0 & 0 & \cdots & 0 & 0\cr
                   \cdot & \cdot & \cdot & \cdots & \cdot & \cdot\cr
                   0 & 0 & 0 & \cdots & 1 & 0\cr
                   0 & 0 & 0 & \cdots & 0 & 1\cr
                   0 & 0 & 0 & \cdots & 0 & 0\cr}\right)_{d\times d},  \lb{3.4}\ee
and $A^T$ is the transpose of $A$.

Consider the following initial value problem of the linear Hamiltonian system
\bea
&& \dot{\ga}(t) = JQ\ga(t) \qquad \forall\;t\in [0,1],   \lb{3.5}\\
&& \ga(0) = I. \lb{3.6}\eea

Here and below, according to Examples 2) on p.5 of \cite{AG}, the standard symplectic matrix on the coordinate
space $\R^{2n}$ is given by $J = \left(\matrix{ 0 & -I_n\cr
                                                I_n & 0\cr}\right)$ where $I_n$ is the $n\times n$
identity matrix on $\R^n$. According to the Corollary in p.5 of \cite{AG}, the Darboux coordinates is
written as $(p_1,\ldots,p_d,q_1,\ldots,q_d) \in \R^{2n}$. Note that these choices coincides with those
chosen in Long's book \cite{Lo02} of 2002 which are related to the orientation in the discussions here,
and are important in the computations blow.

Then from (\ref{3.4})-(\ref{3.6}) we obtain
\be  \ga(t) = \left(\matrix{ B^{\ast}(t) & O\cr
                            O & B(t)\cr}\right)_{2d\times 2d}, \qquad \forall\;t\in [0,1], \lb{3.7}\ee
where
$$ B(t) = \left(\matrix{1 & t & t^2/2 & t^3/(3!) & \cdots & \cdots & t^{d-2}/((d-2)!) & t^{d-1}/((d-1)!) \cr
                        0 & 1 & t & t^2/2 & \cdots & \cdots & t^{d-3}/((d-3)!) & t^{d-2}/((d-2)!) \cr
                        0 & 0 & 1 & t & \cdots & \cdots & t^{d-4}/((d-4)!) & t^{d-3}/((d-3)!) \cr
                        0 & 0 & 0 & 1 & \cdots & \cdots & t^{d-5}/((d-5)!) & t^{d-4}/((d-4)!) \cr
                        \cdot & \cdot & \cdot & \cdots & \cdots & \cdot & \cdot \cr
                        0 & 0 & 0 & 0 & \cdots & \cdots & t & t^2/2 \cr
                        0 & 0 & 0 & 0 & \cdots & \cdots & 1 & t \cr
                        0 & 0 & 0 & 0 & \cdots & \cdots & 0 & 1\cr}\right)_{d\times d},   $$
and
{\footnotesize $$ B^{\ast}(t)=\left(\matrix{
  1 & 0 & 0 & 0 & \cdots & \cdots & 0 & 0\cr
  (-t) & 1 & 0 & 0 & \cdots & \cdots & 0 & 0\cr
  (-t)^2/2 & (-t) & 1 & 0 & \cdots & \cdots & 0 & 0\cr
  (-t)^3/3! & (-t)^2/2 & (-t) & 1& \cdots & \cdots & 0 & 0\cr
  \cdot & \cdot & \cdot & \cdot & \cdots & \cdots & \cdot & \cdot\cr
  (-t)^{d-3}/((d-3)!) & (-t)^{d-4}/((d-4)!) & (-t)^{d-5}/((d-5)!) & (-t)^{d-6}/((d-6)!) & \cdots & \cdots & 0 & 0\cr
  (-t)^{d-2}/((d-2)!) & (-t)^{d-3}/((d-3)!) & (-t)^{d-4}/((d-4)!) & (-t)^{d-5}/((d-5)!) & \cdots & \cdots & 1 & 0\cr
  (-t)^{d-1}/((d-1)!) & (-t)^{d-2}/((d-2)!) & (-t)^{d-3}/((d-3)!) & (-t)^{d-4}/((d-4)!) & \cdots & \cdots & (-t) & 1\cr}\right)_{d\times d}. $$}
Here we have
\be  \det(\lm I_{2d} - \ga(t)) = (\lm - 1)^{2d}, \qquad \forall\;t\in [0,1].  \lb{3.8}\ee
Thus $\sg(\ga(t)) = \{1,1,\ldots,1\}$, i.e., $1$ is an eigenvalue of $\ga(t)$ with algebraic multiplicity $2d$ for every $t\in [0,1]$.

By the study in Section 1.4 on pp.16-24 and Theorem 1.8.10 on p.41 of Long's book \cite{Lo02}, the
matrix $\ga(1)$ can be connected to a decomposition of basic normal forms in $\Om_1^0(\ga(1))$,
where the homotopy set $\Om_1(M)$ respect to the eigenvalue $1$ is defined by
\be  \Om_1(M) = \{N\in \Sp(2d)\;|\; \nu_1(N)=\nu_1(M)\},   \lb{3.9}\ee
and $\Om_1^0(M)$ is the path connected component of $\Om_1(M)$ containing $M$ defined in \cite{Lo99b}
and can be find in Definition 1.8.5 on p.38 of \cite{Lo02}. Specially the above matrix
$\ga(t)$ with $t\in [0,1]$ is in the form of the matrix $N_d(1,0)$ defined in (1.4.1)-(1.4.5) with
$b=0$ on pp.17-18 in \cite{Lo02} which is written as
$$    N_d(1,0)= \left(\matrix{A_d(1) & O\cr
                              O & C_d(1)\cr}\right)_{2d\times 2d}.   $$
Thus by Theorem 1.4.1 on p.18 of \cite{Lo02}, we obtain
\be  \nu(\ga) = \dim\ker (\ga(1)-I) = 2.  \lb{3.10}\ee
where and below we write $\nu(\ga)=\nu_1(\ga)$ for notational simplicity, when the eigenvalue $1$
in consideration is clear enough.

Using the basic normal forms (defined in \cite{Lo99a}, \cite{Lo00}) as in (1.8.15)-(1.8.17) on p.41 in
\cite{Lo02}, by Theorem 1.8.10 on p.41 of \cite{Lo02} again, it yields
\be      I_2 \dm H_0 \in \Om_1^0(\ga(1)),  \lb{3.11}\ee
where $H_0\in \Sp(2(d-1))$ satisfies $\sg(H_0)\cap\U = \emptyset$.

In this case, authors of \cite{CGG24} defined $b_0(Q)$ to be the number of $Q_0$ in the decomposition
of $Q$. When $Q=Q_0$ given by (\ref{3.1}) they defined
\be    b_0(Q_0)=1.   \lb{3.12}\ee
Consequently by the definition (\ref{3.3}), it yields
\be   \beta_{\pm}(Q_0) = b_0(Q_0) = 1. \lb{3.13}\ee

Using the splitting numbers $S^{\pm}_M(\lm)$ defined in \cite{Lo99b} (cf. Definition 9.1.4 on p.191
of \cite{Lo02}), by the List 12 on p.198 of \cite{Lo02} we have
\be   S_{\ga(1)}^{\pm}(1) = S_{I_2}^{\pm}(1) = 1. \lb{3.14}\ee
Thus, because both $S_M^{\pm}(\lm)$ and $\beta_{\pm}(M)$ are additive in the decomposition of $M$
into sums of symplectic matrices, in this case we have
\be   S_{\ga(1)}^{\pm}(1) = \beta_{\pm}(\ga). \lb{3.15}\ee

\medskip

{\bf Step 1.3.} {\it On $b_{\ast}(\Phi)$ with $\ast = \pm$.}

Next we study the quadratic form $Q_{\pm}$ on $\R^{2d}$ with an odd integer $d\ge 1$.

Because the quadratic forms $Q_0$ and $Q_{\pm}$ are operating on different subspaces of $\R^{2d}$, we
can study them respectively, and then add them together to get the total quadratic form finally.

{\bf Subcase 1.3.1.} {\it The subcase $Q_0=0$.}

In this subcase, by the definition in \cite{AG} and \cite{CGG24}, the quadratic form $Q_{\pm}$ is defined
on $\R^2$ as $Q_{\ep}$ with $\ep=1$ or $-1$ and the coordinates on $\R^2$ is given by $(p_d,q_d)$.
In this case, from
$$  Q_{\ep}\left(\matrix{p_d\cr
                         q_d\cr}\right)\cdot \left(\matrix{p_d\cr
                                                           q_d\cr}\right)  = \ep \frac{p_d^2}{2}, $$
we obtain
\bea  Q_{\ep} = \left(\matrix{\ep & 0\cr
                             0   & 0\cr}\right).   \nn\eea
Then we consider the following initial value problem of the linear Hamiltonian system and we denote
its fundamental solution by $\ga_{\ep}(t)$ with $\ep=1$ or $-1$,
\bea
&& \dot{\ga}_{\ep}(t) = JQ_{\ep}\ga_{\ep}(t), \qquad \forall\;t\in [0,1], \lb{3.16}\\
&&  \ga_{\ep}(0) = I_2.    \lb{3.17}\eea

Then solving this system directly, we obtain
\be  \ga_{\ep}(t) =  \left(\matrix{1 & 0\cr
                             \ep t  & 1\cr}\right), \qquad \forall\;t\in [0,1]. \lb{3.18}\ee
Then
\be  \ga_{\ep}(1) =  \left(\matrix{1 & 0\cr
                                    \ep & 1\cr}\right) \in \Om_1^0(N_1(1,-\ep)),  \lb{3.19}\ee
with $N_1(1,-\ep) =  \left(\matrix{1 & -\ep\cr
                              0 & 1\cr}\right)$. Therefore by List 9.1.12 on p.198 of \cite{Lo02},
we obtain
\bea
&& \nu(\ga_{-1}) = \nu(N_1(1,1)) = 1,\qquad S_{\ga_{-1}(1)}^{\pm}(1) = S_{N_1(1,1)}^{\pm}(1) = 1, \lb{3.20}\\
&& \nu(\ga_{1}) = \nu(N_1(1,-1)) = 1,\qquad S_{\ga_{1}(1)}^{\pm}(1) = S_{N_1(1,-1)}^{\pm}(1) = 0. \lb{3.21}\eea

Specially because $\ga_1(1)\in \Om_0(N_1(1,-1))$ and $\ga_{-1}(1)\in \Om_0(N_1(1,1))$ hold, it
yields
\be \bb_+(\ga_1) = 1, \quad\bb_-(\ga_1) = 0, \quad {\rm and}\quad S_{\ga_1(1)}^{\pm}(1) = S_{N_1(1,-1)}^{\pm}(1) = 0, \lb{3.22}\ee
hold simultaneously, and
\be \bb_+(\ga_{-1}) = 0, \quad\bb_-(\ga_{-1}) = 1, \quad {\rm and}\quad S_{\ga_{-1}(1)}^{\pm}(1) = S_{N_1(1,1)}^{\pm}(1) = 1, \lb{3.23}\ee
hold simultaneously too. Thus, in this subcase we have always
\be   S_{\ga_\epsilon(1)}^{-}(1) = \beta_-(\ga_{\epsilon}) \qquad {\rm for}\quad \ep = \pm 1. \lb{3.24}\ee

\medskip

{\bf Subcase 1.3.2.} {\it The subcase $Q_0 \not= 0$.}

In this subcase, for $\ep = \pm 1$, by the definition in \cite{AG} and \cite{CGG24} the quadratic form
$Q_{\ep}$ is given by $Q_{\ep} = \ep(Q_0 + p_d^2/2)$. Let $\ga_{\ep}(t)$ be the fundamental solution of
the initial value problem of the linear Hamiltonian system (\ref{3.16})-(\ref{3.17}).

Then considering the above study on $Q_0$ on $\R^{2d}$ and $Q_{\pm}$ on $\R^2$, specially (\ref{3.3}) and
(\ref{3.10}), by the study in Sections 1.8 and 1.9 on pp.36-47 in \cite{Lo02}, the basic normal form
decomposition of $\ga_{\pm 1}$ can be written as follows,
\be   I_2 \dm N_1(1,\mp 1) \dm H_{\mp 1} \in \Om^0(\ga_{\pm 1}(1)),  \lb{3.25}\ee
where $H_{\mp 1}\in \Sp(2(d-2))$ satisfies $\sg(H_{\mp 1})\cap \U = \emptyset$.

Because the invariants in (\ref{3.3}) and the splitting numbers $S_M^{\pm}(\lm)$ in \cite{Lo99b} and
Definition 9.1.4 on p.191 of \cite{Lo02} are symplectic additive, by (\ref{3.25}) in the subcase
1.3.2 we obtain
\be   S_{\ga_{\ep}(1)}^{-}(1) = \beta_-(\ga_{\ep}). \lb{3.26}\ee

\medskip

Summarising the above studies, using the symplectic additivity of related invariants, we obtain

{\bf Proposition 1.} {\it For any $\ga\in C([0,1],\Sp(2n))$, by the above proofs, the following
holds always, }
\be  \bb_{-}(\ga) = S_{\ga(1)}^-(1).  \lb{3.27}\ee

\medskip

\setcounter{equation}{0}\setcounter{figure}{0}
\section{The (IR1) in Theorem 3.3 of \cite{CGG24} is a direct consequence of Remark 3.6 of \cite{DLW16}} %Section 4

In Section 3 of \cite{CGG24}, as usual the authors denote the set of all non-degenerate
symplectic matrices in $\Sp(2m)$ with $m\in\N$ by $\Sp^{\ast}(2m)$. Then for a symplectic path $\Phi$
with end point $\Phi(1)\in \Sp^{\ast}(2m)$, the end point matrix $\Phi(1)$ can be connected to $-I$ in
$\Sp^{\ast}(2m)$ by a path $\Psi$ lying entirely in $\Sp^{\ast}(2m)$. Let $\Phi'$ be the concatenating
path of $\Phi$ and $\Psi$. Then as usual the Conley-Zehnder index of $\Phi$ is defined by
\be   \mu(\Phi) = \hat{\mu}(\Phi'),    \lb{4.1}\ee
where $\hat{\mu}$ is defined to be the unique quasimorphism $\hat{\mu}: \widetilde{\Sp}(2m) \to \R$
which is continuous, homogeneous and satisfies the normalized condition:
$$   \hat{\mu}(\Phi_0) = 2 \qquad {\rm for}\quad \Phi_0(t) = \exp(2\pi\sqrt{-1} t)\oplus I_{2m-2}. $$
and $\widetilde{\Sp}(2m)$ is the universal covering of $\Sp(2m)$.

Then the authors of \cite{CGG24} defined two indices for $\Phi\in \widetilde{\Sp}(2m)$ by
\be  \mu_+(\Phi) := \lim\sup_{\bar{\Phi}\to \Phi} \mu(\bar{\Phi}), \quad
    \mu_-(\Phi) := \lim\inf_{\bar{\Phi}\to \Phi} \mu(\bar{\Phi}),   \lb{4.2}\ee
where in both cases the limit is taken over $\bar{\Phi}\in \widetilde{\Sp^*}(2m)$ converging to $\Phi$.

Our next result shows that the relations of the property (IR1") with the indices defined in \cite{Lo02} as
well as \cite{DLW16} are as follows,

{\bf Proposition 2.} {\it (i) The assertion (IR1"-1) in (IR1) are direct consequence of (3.42)-(3.43) in
\cite{DLW16}.

(ii) The assertion (IR1"-2) in (IR1) was already proved in (1) of Theorem 1.1 of \cite{LL00} (cf. $1^0$ of
Theorem 10.1.2 on p.213 of \cite{Lo02}). }

{\bf Proof.} Note first that using definitions of $\mu_{\pm}(\Phi)$ in (\ref{4.2}), by Corollary 6.1.12
and Definition 6.1.13 on p.144 of \cite{Lo02}, it yields immediately
\bea
&&  \mu_-(\Phi) = i(\Phi), \lb{4.3}\\
&&  \mu_+(\Phi) = i(\Phi) + \nu(\Phi). \lb{4.4}\eea
Note that such identities were first introduced in \cite{Lo90} of 1990 and rigorously proved in
\cite{Lo97} of 1997 by Long. Then they yield
\be    \hat{\mu}(\Phi) = \hat{i}(\Phi), \qquad \forall\;\Phi\in \Sp(2m). \lb{4.5}\ee

To prove (IR1"-1), letting $k$ and $d$ in (IR1"-1) be given by $k=2m_i$ and $d=2N$. Note first that
using notations in \cite{Lo02} and \cite{DLW16}, the left hand side of (IR1"-1) is given by
$$  |\hat{\mu}(\Phi_i^{k}) - d| = |2m_i\hat{i}(\Phi_i) - 2N|. $$
where $2m_i\hat{i}(\Phi_i)= \hat{i}(\Phi_i^{2m_i})$ holds. Then by (3.42) and (3.43) of \cite{DLW16},
for any given $\eta>0$, we obtain
\bea  |2m_i\hat{i}(\Phi_i) - 2N|
&=& \left|2([\frac{N}{\bar{M}\hat{i}(\Phi_i)}] + \chi_i)\bar{M}\hat{i}(\Phi_i)
                 - 2 \frac{N}{\bar{M}\hat{i}(\Phi_i)} \bar{M}\hat{i}(\Phi_i)\right| \nn\\
&=& \left|-2\bar{M}\hat{i}(\Phi_i)\left(\{\frac{N}{\bar{M}\hat{i}(\Phi_i)}\} - \chi_i\right)\right|   \nn\\
&<& 2 \bar{M}\hat{i}(\Phi_i)\ep  < \eta, \nn\eea
when we take $\ep>0$ small enough. Therefore (IR1"-1) is proved.

Next using notations in \cite{Lo02} and \cite{DLW16}, the integer $d$ in (IR1"-2) is given by $d = 2N$ too.
Therefore (IR1"-2) coincides precisely with (1) of Theorem 1.1 of \cite{LL00} (cf. $1^0$ of
Theorem 10.1.2 on p.213 of \cite{Lo02}).

The proof is complete. \hfill\hb

\medskip

\setcounter{equation}{0}\setcounter{figure}{0}
\section{The coincidence of (IR2) and (IR3) in Theorem 3.3 of \cite{CGG24} and (3.20)-(3.22) in
Theorem 3.5 of \cite{DLW16}} %Section 5

Now in the studies of the next two propositions, we prove that (IR2) and (IR3) of \cite{CGG24}
coincide completely with (3.20)-(3.22) of Theorem 3.5 in \cite{DLW16}. In these proofs we use the
notations in \cite{DLW16} and \cite{Lo02}.

\medskip

{\bf Proposition 3.} {\it The assertion (IR2") and then (IR2) in \cite{CGG24} coincides precisely with
(3.20)-(3.21) in \cite{DLW16}. }

{\bf Proof.} Note first by Corollary 6.1.12 and Definition 6.1.13 on p.144 of \cite{Lo02}, it yields
immediately (\ref{4.3})-(\ref{4.5}). Therefore (IR2") and then (IR2) in \cite{CGG24} coincides precisely
with (3.20)-(3.21) of \cite{DLW16}. \hfill\hb

\medskip

Next we have

\medskip

{\bf Proposition 4.} {\it The assertion (IR3") and then (IR3) in \cite{CGG24} coincides precisely with
(3.20)-(3.22) in \cite{DLW16}. }

\medskip

{\bf Proof.} In order to use uniformly the notations in \cite{DLW16} in the following proof, we have changed
the notations $\R^{2m}$, $\Sp(2m)$, $\ell$ and $\Phi$ in (IR3") of \cite{CGG24} to $\R^{2n}$, $\Sp(2n)$, $m$ and
$\ga$ as those in (3.20) and (3.22) of \cite{DLW16} respectively,

Note first that by the requirement of Theorem 3.3 of \cite{CGG24} on the choice of $k$, the degree $D$ of
$\ga_k$ is a factor of $k$, and by the requirement (3.30) on $m_k$ of Theorem 3.5 of \cite{DLW16},
there holds $\frac{m_k\th}{\pi}\in\Z$ whenever $\frac{\th}{\pi}\in \Q\cap (0,2)$ and
$e^{\sqrt{-1}\th}\in \sg(M_k)$. These two requirements are precisely the same.

By (IR3") in \cite{CGG24} and (3.22) in \cite{DLW16}, to prove Proposition 4, it suffices to prove
the following identity for every symplectic path in $C([0,1],\Sp(2n))$,
\be   \bb_+(\ga_k^m) - \bb_-(\ga_k^m) = \nu(\ga_k^{2m_k-m}) - 2(S_{M_k}^+(1) + Q_k(m)), \lb{5.1}\ee
where $M_k=\ga_k(1)$, and $Q_k(m)$ is given by (3.25) of \cite{DLW16}.

By the above identities (\ref{3.2}) and (\ref{3.3}), we obtain
\be   \bb_+(\ga_k^m) +\bb_-(\ga_k^m)=\nu(\ga_k^{m})=\nu(\ga_k^{2m_k-m}),  \lb{5.2}\ee
where the last equality follows from (3.20) of \cite{DLW16}.

Then subtracting (\ref{5.2}) from (\ref{5.1}) yields
\be   \bb_-(\ga_k^m) = S_{M_k}^+(1) + Q_k(m) = S_{M_k}^-(1) + Q_k(m), \lb{5.3}\ee
where to get the last equality we have applied the fact $S_{M_k}^+(1)=S_{M_k}^-(1)$
from Lemma 9.1.6 on p.192 of \cite{Lo02}.

Then by the Bott-type formula for splitting numbers which was first proved in \cite{LLZ02} by Liu,
Long and Zhu (cf. Corollary 9.2.4 on p.201 in \cite{Lo02}), we obtain
\be S_{M_k}^-(1) + Q_k(m) = \sum_{\omega^m=1}S_{M_k}^-(\omega) = S_{M_k^m}^-(1), \lb{5.4}\ee

On the other hand, by Proposition 1 and the fact $M_k^m = \ga_k^m(1)$, we have
$$  \bb_-(\ga_k^m) = S_{M_k^m}^-(1).   $$
Together with (\ref{5.4}), it yields a proof of (\ref{5.3}) and then (\ref{5.1}).

The proof of Proposition 4 is complete. \hfill\hb

\medskip

\bibliographystyle{abbrv}

\emph{}

\end{document}